\documentclass[10pt]{article}

\usepackage{amssymb,latexsym,amsmath,amsthm}
\usepackage{graphicx} % Replaced 'epsfig', as 'graphicx' is more modern.
\usepackage{serre}
\usepackage{url}
\usepackage{authblk}  % Added for author affiliations

\title{On a remark of Serre}

\author[1]{Katie Ahrens}
\author[2]{Jon Grantham}

\affil[1]{Institute for Defense Analyses\\ Center for Communications Research\\ La Jolla, CA}
\affil[2]{Institute for Defense Analyses\\ Center for Computing Sciences\\ Bowie, MD}

\date{}
\newtheorem{proposition}{Proposition}

\begin{document}

\maketitle

\begin{abstract}
Inspired by a remark of Serre, we extend the search for primes $p$ such that the maximum Hasse bound for the number of points on an elliptic curve over $\mathbb{F}_{p^5}$ is not achieved. We then give a list of all $q<10^{70}$ such that the Hasse bound is not achieved over $\mathbb{F}_{q}$. We explore the heuristics for how many such numbers should exist in each case. Finally, we look at similar criteria for genus $2$ and $3$ curves.
\end{abstract}

\section{Introduction}

Hasse \cite{hasse} proved that the maximum number of points on an elliptic curve over the finite field $\mathbb{F}_q$ with $q$ elements is $q+\lfloor{2\sqrt{q}}\rfloor+1$. Weil \cite{weil} generalized this by proving the Riemann hypothesis for curves over finite fields, which yields $N_q(g)\le q+1+2g\sqrt{q}$. Serre sharpened this to $N_q(g)\le q+g\lfloor{2\sqrt{q}}\rfloor+1$. Denote the maximum number of points on any smooth projective curve over $\mathbb{F}_{q}$ of genus $g$ by $N_q(g)$. Let $m=\lfloor 2\sqrt{q} \rfloor$.

If the Hasse-Weil bound is tight for curves of genus $g$ over $\F_q$ (that is, $N_q(g)=q+gm+1$), we say that the curves over $\F_q$ have defect 0. We may shorten this and say that ``$q$ has defect $0$'' when the context is clear. In general, we define the defect of $q$ to be $(q+gm+1)-N_q(g)$. 
 
Consider the genus $1$ case. Serre \cite{serre} notes that results of Deuring \cite{deuring} and Waterhouse \cite{waterhouse} imply a condition for when this maximum
is achieved. Namely, let $q=p^e$. The Hasse maximum fails to be achieved when $p \mid \lfloor{2\sqrt{p^e}}\rfloor$ with $e$ odd, except in the two cases $(p,e)=(2,1)$ and $(3,1)$; in all failure cases, the defect is $1$. Serre observes that this
implies $e\ge 5$ and, in Remark 2.6.5, says that $p=7$ is the only example with $e=5$ for $p<10^7$. He says that a ``naive log log argument'' implies that there should be ``very sparse'' larger primes.

We propose the term {\it Serre prime} for any prime $p$ such that there is no elliptic curve over $\mathbb{F}_{p^5}$ with $p^5+\lfloor{2\sqrt{p^5}}\rfloor+1$ points and the term {\it Deuring-Waterhouse number} for a prime power $q$ such that there is no elliptic curve over $\mathbb{F}_q$ with $q+\lfloor{2\sqrt{q}}\rfloor+1$ points.

In this note, we extend the search for Serre primes up to $10^{16}$ and enumerate all Deuring-Waterhouse numbers up to $10^{70}$. We use
data from the search to evaluate the plausibility of the log log argument.

In Sections \ref{sec:genus2} and \ref{sec:genus3}, we look at similar criteria for genus $2$ and genus $3$ curves, respectively.
The genus $2$ case is similar to the genus $1$ case in that we can completely characterize the prime powers with nonzero defect using some simple arithmetic tests.
In genus $3$, no complete characterization is known, although the concept of ``minimal relative defect'' will allow us to present some partial results.

\section{Search for Serre primes}

We generate the primes up to $10^6$ with a sieve of Eratosthenes. Then, we divide the integers up to $10^{16}$ into congruence classes modulo $510510$. This allows us to parallelize the search. Sieving by the primes up to $10^6$
allows us to remove most composites. We use GMP \cite{gmp} to compute $m$ and see if it is divisible by $p$; if it is, we perform a primality test on $p$. We found no primes $p$ dividing $\lfloor 2\sqrt{p^5}\rfloor$ other than $7$.

\section{The naive log log argument}

	It is natural to ask whether or not we should have expected to find any Serre primes other than $7$. Expanding Serre's log log argument gives a heuristic.

Without any known obstruction, we assume that $m=\lfloor 2\sqrt{p^5}\rfloor$ behaves as if chosen randomly mod $p$. In particular, we assume that $p\mid m$ with probability $1/p$. So the expected number of $p$ with $y<p<x$ and $p\mid m$ is $\sum_{y<p<x} 1/p \approx \log\log(x)-\log\log(y)$, where the approximation follows from Mertens' second theorem.
$\log\log(10^{16})-\log\log(10^{7})\approx 0.83$, so our failure to find another Serre prime was not just due to spectacularly bad luck.

It seems unlikely further computation in the $e=5$ case will allow us to test the heuristic. But the log log argument should apply for all odd $e\ge 5$. So we looked at $7\le e < 107$. We searched $p$ with $10^{2}<p<10^{11}$, and
$\log\log(10^{11})-\log\log(10^{2})\approx 1.7$. There are $50$ odd values of $e$ in this range, so we expect approximately $50\cdot 1.7=85$ primes total. We found $84$.

\section{Enumeration of Deuring-Waterhouse numbers}

These numbers are listed as sequence A364690 of the OEIS \cite{oeis}. We have updated them with our calculations. When searching for $p^e<B$ for $e\ge 7$, we use the same algorithm as
in the previous section.

\begin{proposition}[Genus 1 results]
\label{prop:genus1}
For genus $g=1$, the values of $q<10^{70}$ having positive defect are precisely the following $146$ Deuring--Waterhouse numbers:

$2^{7}$, $2^{11}$, $2^{15}$, $2^{17}$, $2^{19}$, $2^{21}$, $2^{23}$, $2^{27}$, $2^{29}$, $2^{39}$, $2^{41}$, $2^{47}$, $2^{49}$, $2^{55}$, $2^{57}$, $2^{73}$, $2^{75}$, $2^{83}$, $2^{93}$, $2^{95}$, $2^{101}$, $2^{103}$, $2^{107}$, $2^{109}$, $2^{113}$, $2^{115}$, $2^{117}$, $2^{119}$, $2^{123}$, $2^{125}$, $2^{127}$, $2^{131}$, $2^{137}$, $2^{139}$, $2^{143}$, $2^{155}$, $2^{161}$, $2^{163}$, $2^{165}$, $2^{169}$, $2^{171}$, $2^{177}$, $2^{183}$, $2^{185}$, $2^{191}$, $2^{199}$, $2^{203}$, $2^{207}$, $2^{211}$, $2^{213}$, $2^{217}$, $2^{221}$, $2^{225}$, $2^{237}$, $2^{241}$, $2^{243}$, $2^{255}$, $2^{257}$, $2^{259}$, $2^{267}$, $2^{271}$, $2^{277}$, $2^{287}$, $2^{293}$, $2^{295}$, $2^{297}$, $2^{299}$, $3^{7}$, $3^{15}$, $3^{25}$, $3^{37}$, $3^{41}$, $3^{43}$, $3^{47}$, $3^{53}$, $3^{65}$, $3^{69}$, $3^{75}$, $3^{85}$, $3^{87}$, $3^{91}$, $3^{97}$, $3^{103}$, $3^{105}$, $3^{107}$, $3^{109}$, $3^{143}$, $3^{149}$, $3^{151}$, $3^{153}$, $3^{159}$, $3^{177}$, $3^{181}$, $3^{185}$, $3^{191}$, $3^{201}$, $3^{211}$, $3^{213}$, $3^{217}$, $3^{219}$, $3^{223}$, $3^{225}$, $3^{227}$, $3^{229}$, $5^{9}$, $5^{11}$, $5^{15}$, $5^{19}$, $5^{27}$, $5^{37}$, $5^{39}$, $5^{51}$, $5^{69}$, $5^{95}$, $5^{111}$, $5^{113}$, $5^{117}$, $5^{127}$, $5^{145}$, $7^{5}$, $7^{29}$, $7^{49}$, $7^{81}$, $7^{83}$, $7^{91}$, $13^{17}$, $13^{19}$, $13^{29}$, $13^{39}$, $13^{45}$, $13^{51}$, $13^{63}$, $17^{41}$, $23^{35}$, $37^{19}$, $79^{23}$, $109^{27}$, $113^{9}$, $131^{33}$, $239^{9}$, $1249^{19}$, $43783^{9}$, $73877^{11}$, $248461^{11}$, $747343^{7}$, $411484343^{7}$.

All of these have defect $1$; all other $q<10^{70}$ have defect $0$.
\end{proposition}

We ask ourselves a similar question to what we asked in the case of the Serre primes --- is this the sort of pattern we should have been expecting?
Serre provides a short proof that there are an infinite number of base-2 Deuring-Waterhouse numbers by considering the base-2 expansion of $\sqrt{2}$. 
Each $0$ in this expansion corresponds to a power $q$ of $p=2$ such that $p$ divides $\lfloor 2q^{1/2} \rfloor$ (since multiplying $q$ by $p^2$ shifts the base-$p$ expansion of $2q^{1/2}$ by one place, and the floor truncates the rest of the expansion).
In the base $2$ case, it is obvious that $0$ appears in this expansion an infinite number of times, since otherwise $\sqrt{2}$ would be rational.

We can in general consider the base $p$ expansion of $2\sqrt{p}$ and attempt to make the same argument. But we cannot say much more than that, since it is not obvious 0 appears an infinite number of times in these other expansions.
Although it is interesting to note that Borel's 1950 conjecture \cite{borel} that every irrational algebraic number is normal would imply that every prime $p$ will eventually appear as a base in the list of Deuring-Waterhouse numbers an infinite number of times.

\section{Genus 2}
\label{sec:genus2}

\subsection{Background}

In genus $1$, a simple divisibility check sufficed to see if $q$ has positive defect. In genus $2$, Serre tells us that there are now two ways for $q$ to have positive defect. The first is through a divisibility condition, and the second is if $q$ can be represented as a certain polynomial. We recall Serre's results.

First, everything is known in the case of $q$ a square, by Theorem 3.2.3 of \cite{serre}: the defect is $0$, unless $q=9$ in which case the defect is $2$, or $q=4$ in which case the defect is $3$. 

We have the following theorem (\cite{serre} Definition 3.2.1 and Theorem 3.2.3).

\begin{theorem}
Let $q=p^e$ for $p$ a prime, and suppose $q\notin\{4,9\}$. Let $\{k\}$ denote the fractional part of a number $k$. We say $q$ is special if it is not a square, and it satisfies one of the following two conditions:
\begin{enumerate}
\item $p$ divides $m:=\lfloor 2\sqrt{q}\rfloor$
\item $q$ is represented by one of the polynomials $x^2+1, x^2+x+1, x^2+x+2$.
\end{enumerate}

Then:
\begin{itemize}
\item $N_q(2)=q+2m+1$ if $q$ is nonspecial.
\item $N_q(2) = q+2m$ if $q$ is special and $\{2q^{1/2}\} > \frac{\sqrt{5}-1}{2}$.
\item $N_q(2) = q+2m-1$ if $q$ is special and $\{2q^{1/2}\} < \frac{\sqrt{5}-1}{2}$.
\end{itemize}
\end{theorem}

In Theorem 3.2.6, Serre proves a remarkably powerful result about the case where $q=p^e$ with $e\ge 2$: the only $q$ that are special by condition $(2)$ of the definition are $7^3, 2^3, 2^5, 2^{13}$.

So we are left with two avenues of computational exploration of positive defect curves in genus $2$. 

\begin{enumerate}
\item We can find $q=p^e$, $e\ge 2$ that are special by condition (1).
\item We can find $q=p$ that are represented by $x^2+1$ or $x^2+x+1$ (note $x^2+x+2$ is always even).
\end{enumerate}
The prime powers satisfying condition (1) are exactly the Deuring-Waterhouse numbers, which we already computed in the $g=1$ case.

\begin{proposition}[Genus 2 results]
\label{prop:genus2}
For genus $g=2$:
\begin{enumerate}
\item There are $61$ non-prime $q<10^{70}$ with defect $1$, and $85$ non-prime $q<10^{70}$ with defect $2$. These are the Deuring-Waterhouse numbers, listed in Proposition~\ref{prop:genus1}.
\item $q=2$ has defect $1$.
\item The primes $q>2$ of the form $q=x^2+1$ or $q=x^2+x+1$ all have defect $2$, while $q=2=1^2+1$ has defect $1$. The number of primes $q<10^n$ of each of these two forms, for $n=1,\dots,22$, is given in the corresponding columns of Table~\ref{tab:primecounts}; these counts include $q=2$ in the $x^2+1$ column.
\end{enumerate}
\end{proposition}

In the $q=p$ case, note first that the results of the triple sieve of the second author and Graves in \cite{granthamgraves2025} give us an exhaustive list of primes $p<6.25\times 10^{28}$ of the form $x^2+1$. 

We have adapted this same code to the $x^2+x+1$ case.  Note that for a prime $p>3$, 
$$x^2+x+1 \equiv 0 \mod p \iff x^3 \equiv 1 \mod p \text{ and } x\not\equiv 1 \mod p.$$

So, in this case, if we want to identify all primes $p<B$ of the form $x^2+x+1$, we start by generating a list of primes up to $B^{1/4}$ using
the sieve of Eratosthenes. We use that list to sieve up to $B^{1/2}$ to find all primes that are $1\bmod 3$. For the third sieve, we need
the corresponding roots of $x^2+x+1$ for the primes on our second list. We find these by choosing $g\bmod p$ until we find one where $g^{(p-1)/3}\not\equiv 1$; that will necessarily be a root of $x^2+x+1$. Finally, our third sieve looks for values of $x$ such that $x^2+x+1$ is not
divisible by any prime less than $B^{1/2}$. Here, we need to use the roots of $x^2+x+1$ modulo the primes that are $1\bmod 3$.

Table~\ref{tab:primecounts} gives counts of primes represented by the relevant quadratic polynomials. Note that Poletti \cite{poletti} did the first computations of primes of the form $x^2+x+1$. We are not aware of any further published computations in the near-century since, although Bateman and Stemmler \cite{batestem} mention Horn's computations in the case where $x$ is prime.

\begin{table}[h!]
\centering
\begin{tabular}{|| r | r | r | r | r || }
\hline
$B$ & $x^2+x+1$ & $x^2+1$ & $x^2+2$ & $x^2+x+3$ \\
\hline
$10^{1}$ & $2$ & $2$ & $2$ & $2$ \\
$10^{2}$ & $6$ & $4$ & $4$ & $4$ \\
$10^{3}$ & $14$ & $10$ & $6$ & $8$ \\
$10^{4}$ & $32$ & $19$ & $12$ & $14$\\
$10^{5}$ & $76$ & $51$ & $28$ & $37$ \\
$10^{6}$ & $189$ & $112$ & $69$ & $93$ \\
$10^{7}$ & $520$ & $316$ & $162$ & $245$ \\
$10^{8}$ & $1410$ & $841$ & $447$ & $629$\\
$10^{9}$ & $3825$ & $2378$ & $1237$ & $1708$ \\
$10^{10}$ & $10751$ & $6656$ & $3423$ & $4899$ \\
$10^{11}$ & $30580$ & $18822$ & $9777$ & $13862$ \\
$10^{12}$ & $88118$ & $54110$ & $27869$ & $40037$ \\
$10^{13}$ & $255408$ & $156081$ & $81143$ & $116097$ \\
$10^{14}$ & $745582$ & $456362$ & $236985$ & $339824$ \\
$10^{15}$ & $2188876$ & $1339875$ & $695999$ & $997259$ \\
$10^{16}$ & $6456835$ & $3954181$ & $2054022$ & $2941031$ \\
$10^{17}$ & $19145042$ & $11726896$ & $6090628$ & $8714892$ \\
$10^{18}$ & $56988601$ & $34900213$ & $18127693$ & $25941031$  \\
$10^{19}$ & $170218020$ & $104248948$ & $54151650$ & $77483932$  \\
$10^{20}$ & $510007598$ & $312357934$ & not computed & not computed \\
$10^{21}$ & $1532303073$ & $938457801$ & not computed & not computed \\
$10^{22}$ & $4615215645$ & $2826683630$ & not computed & not computed \\
\hline
\end{tabular}
\caption{Number of primes $q<B$ represented by the given polynomials, where $x$ ranges over all nonnegative integers. In particular $q=2$ is included, appearing as $1^2+1$ in the $x^2+1$ column and as $0^2+2$ in the $x^2+2$ column.}
\label{tab:primecounts}
\end{table}

The column for $x^2+1$ is from \cite{granthamgraves2025}. Note that there are more primes of the form $x^2+x+1$ than of the form $x^2+1$. The Bateman-Horn conjecture \cite{BH} predicts that the ratio of these two quantities should converge to $\approx 1.626$. At $10^{22}$, it is $1.6327\dots$.

In contrast to the Deuring-Waterhouse numbers, where the prime powers behaved as if uniformly distributed mod $1$, the only prime $q$ of the form $x^2+1$ or $x^2+x+1$ which has $\{2q^{1/2}\} > \frac{\sqrt{5}-1}{2}$ is $2$.
To see why, note that for $q=x^2+1$ with $x\ge 2$ we have $x^2 < q < x^2+\frac{1}{2}x+\frac{1}{16}$, so $2x < 2q^{1/2} < 2x+\frac{1}{2}$, and thus $\{2q^{1/2}\}<\frac{1}{2}$. Likewise, for $q=x^2+x+1$ with $x\ge 1$ we have $x^2+x+\frac{1}{4} < q < x^2+\frac{3}{2}x+\frac{9}{16}$, so $2x+1 < 2q^{1/2} < 2x+\frac{3}{2}$, and thus $\{2q^{1/2}\}<\frac{1}{2}$. In both cases (unless $x=1$ and $q=x^2+1=2$), the fractional part $\{2q^{1/2}\}$ is less than $\frac{1}{2} < \frac{\sqrt{5}-1}{2}$. Thus, for all $q>2$ of either form, the defect is $2$.

\section{Genus 3}
\label{sec:genus3}

The genus $3$ case is less clean: no simple arithmetic tests are known to compute the defect exactly. However, we can still bound it using Serre's notion of ``minimum relative defect.''

From Section 4.4.2 of Serre:

\begin{definition}
The minimal relative defect $a$ is the smallest integer $a$ such that there exists a genus $3$ curve $X/\F_q$ with $|\#X(\F_q) -q-1| = 3m-a$. In particular, the defect is greater than or equal to the minimal relative defect.
\end{definition}

\subsection{Arithmetic setup}

Let $q=p^e$. Lauter's work \cite{lauter2002}, as summarized in Theorem 4.4.3 of \cite{serre}, gives arithmetic tests to compute $a$. If $e$ is even, the minimal relative defect is 0 unless $q=4$ or $16$, in which case $a=3$. So in our computations we focus on the non-square case (that is, odd $e$).

\begin{proposition}[Odd-exponent criterion from \cite{lauter2002}]
\label{prop:genus3odd}
If $e$ is odd, the minimal relative defect $a$ is $0$ unless one of the following holds:
\begin{enumerate}
\item $q=x^2 +r $ with $x$ an integer, $r=1$ or $2$, $r\le x$, for which $a=2$.
\item $q=x^2 + x +r$ with $x$ an integer, $r=1$ or $3$, $r\le x$, for which $a=3$.
\item $p \mid m$, for which $a=3$ unless $\{2q^{1/2}\} \ge 1-4\cos^2(3\pi/7)$, in which case $a=2$.
\end{enumerate}
\end{proposition}

\subsection{Search strategy}

We first treat the proper prime powers, for which a complete classification is available. The following may be viewed as a genus~3 analogue of Theorem~3.2.6 of \cite{serre}.

\begin{proposition}[Proper prime powers]
\label{prop:genus3powers}
Let $q=p^e$ be a proper prime power with $e\ge 3$ odd. If $q$ satisfies condition (1) or (2) of Proposition~\ref{prop:genus3odd}, then $q$ is one of
\[
3^3 = 27 = 5^2+2, \qquad 3^5 = 243 = 15^2+15+3, \qquad 7^3 = 343 = 18^2+18+1.
\]
The first satisfies condition (1) with $a=2$; the latter two satisfy condition (2) with $a=3$.
\end{proposition}

\begin{proof}
A proper prime power $q=p^e$ with $e\ge 3$ is, in particular, a perfect power $y^k$ with $y=p$ and $k=e\ge 3$. We treat in turn the four polynomials appearing in conditions (1) and (2) of Proposition~\ref{prop:genus3odd}.

For $q=x^2+1$ (condition (1) with $r=1$), the result of Lebesgue \cite{Lebesgue} shows that $y^k=x^2+1$ has no solutions with $k\ge 2$, so no proper prime power arises.

For $q=x^2+2$ (condition (1) with $r=2$), a theorem of Ljunggren, later rediscovered by Nagell (see \cite{cohn}), shows that the only solutions of $y^k=x^2+2$ with $k\ge 2$ have $y^k=27$. Since $27=3^3=5^2+2$ with $r=2\le 5=x$, this gives $q=3^3$.

For $q=x^2+x+1=\frac{x^3-1}{x-1}$ (condition (2) with $r=1$), the Nagell--Ljunggren result, in the form given by Bugeaud and Mih\u{a}ilescu \cite[Theorem~NL]{bugeaudmih} (taking $n=3$), shows that the only solutions of $y^k=\frac{x^3-1}{x-1}$ with $k\ge 2$ have $y^k=343$. Since $343=7^3=18^2+18+1$ with $r=1\le 18=x$, this gives $q=7^3$.

For $q=x^2+x+3$ (condition (2) with $r=3$), suppose $y^k=x^2+x+3$ with $k\ge 3$. Multiplying by $4$ gives $4y^k=(2x+1)^2+11$, so that $X^2+11=4y^k$ with $X=2x+1$. If a prime $\ell$ divided both $X$ and $y$, then $\ell\ne 2$ (as $X$ is odd) and $\ell^2\mid 4y^k-X^2=11$, which is impossible; hence $\gcd(X,y)=1$. By Luca, Tengely, and Togb\'e \cite[Theorem~1.1]{lucatt}, $y^k=243$. Since $243=3^5=15^2+15+3$ with $r=3\le 15=x$, this gives $q=3^5$.

These cases exhaust conditions (1) and (2), so the only such $q$ are $3^3$, $3^5$, and $7^3$.
\end{proof}

The value $7^3$ is consistent with the genus~2 prime power result of Serre presented in Section~\ref{sec:genus2}. It remains to count the primes of each relevant form; those of the form $x^2+1$ and $x^2+x+1$ were enumerated in Section~\ref{sec:genus2}.

Because $x^2+2$ and $x^2+x+3$ are not cyclotomic polynomials, we cannot use the root-finding technique to go along with our second sieve from Section~\ref{sec:genus2}, so we perform the computations up to a smaller bound. The resulting prime counts are included in Table~\ref{tab:primecounts}.

\begin{proposition}[Genus 3 results]
\label{prop:genus3}
For genus $g=3$ and odd exponent $e$:
\begin{enumerate}
\item Of the $146$ Deuring-Waterhouse numbers $q<10^{70}$, $26$ have minimum relative defect $2$ and $120$ have minimum relative defect $3$.
\item Among primes, Table~\ref{tab:primecounts} gives the counts of $p<10^n$ represented by $x^2+2$ and by $x^2+x+3$. Those represented with $x\ge r$, that is, all but the primes $2,3$ (for $x^2+2$) and $3,5$ (for $x^2+x+3$), satisfy condition (1), respectively (2), of Proposition~\ref{prop:genus3odd} and so have minimum relative defect $2$, respectively $3$.
\end{enumerate}
\end{proposition}

The known Deuring-Waterhouse numbers with minimum relative defect $2$ can be compared to the expected $28.9$ if they were distributed 
randomly.

Because $1-4\cos^2(3\pi/7)$ is greater than $\frac{\sqrt{5}-1}{2}$ (they are approximately $0.80$ and $0.62$, respectively), a 
Deuring-Waterhouse number with minimum relative defect $2$ in genus $3$ necessarily has defect $1$ in genus $2$, but not vice versa. In fact,
it is easy to deduce from our results that there are $35$ such numbers on our list with defect $1$ in genus $2$ and minimum relative defect
$3$ in genus $3$. The interested reader who wishes to find out the classification of a particular Deuring-Waterhouse number need only compute
the fractional part of $2q^{1/2}$.

\section{A connection to Diophantine equations}
\label{sec:diophantine}

If $p$ is a Serre prime, then $p$ divides $\lfloor 2p^{5/2}\rfloor$. So there is an integer $a$ such that
\[
pa < 2p^{5/2} < pa+1.
\]
Dividing by $p$ and squaring gives
\[
a^2 < 4p^3 < \left(a+\frac1p\right)^2.
\]
Thus $(x,y)=(p,a)$ is an integer solution to 
%Hence, for $x=p$ (prime) and
%\[
%k:=4x^3-z^2,
%\]
%we obtain an integer solution to
\begin{equation}
\label{eq:k_bound}
y^2=4x^3-k
\end{equation}
with
\[
0<k<\frac{2a}{p} + \frac{1}{p^2} < 4\sqrt{p}+\frac{1}{p^2}.
\]

We note equation (\ref{eq:k_bound}) is isomorphic to an elliptic curve in short Weierstrass form with a zero coefficient of the linear term.
If $(p,a)$ is a solution with $a$ even, then $(p,a/2)$ is an integer solution to the equation
\begin{equation*}
y^2=x^3-k/4
\end{equation*}
with $|k/4| < \sqrt{p}+\frac{1}{4p^2}$.

An elliptic curve of the form $y^2=x^3+n, n\neq 0$ is known as a Mordell curve. 
Thus we can reframe the search for Serre primes as a search for integer points $(p,a)$ on the Mordell curve, with the added restriction that the first coordinate be prime and with a restriction on the magnitude of the constant term of the equation relative to $p$.
We are not aware of any results about the number of solutions to a Mordell equation with either of these two restrictions.

\section{Future work}

The tests for genus $1$, $2$, and $3$ all have the same flavor: a divisibility condition and possibly a representation by polynomials, which get slightly more complicated as the genus increases. 
Is there a generalization of such criteria to genus $4$, and maybe beyond?

\section{Acknowledgements}

We thank the participants of the West Coast Number Theory 2024 Conference for helpful suggestions. We especially thank Anay Aggarwal for noting that Serre's base-2 argument generalizes to base-$p$ under the normality assumption.

We thank Jean-Pierre Serre for clarifying the proper attribution of the genus~3 criterion to Lauter.

We also thank the anonymous referee for pointing out the connection between Serre primes and Diophantine equations in Section~\ref{sec:diophantine}, and for suggesting the complete classification of proper prime powers in Proposition~\ref{prop:genus3powers}.

\bibliographystyle{alpha}
\bibliography{serre}

\end{document}